\documentclass[english]{article}
\usepackage[LGR,T1]{fontenc}
\usepackage[latin9]{inputenc}
\usepackage[a4paper]{geometry}
\usepackage{textcomp}
\usepackage{amsmath}
\usepackage{amsthm}
\usepackage{amssymb}
\usepackage{graphicx}
\usepackage{setspace}
\usepackage[authoryear]{natbib}

\makeatletter

\DeclareRobustCommand{\greektext}{%
  \fontencoding{LGR}\selectfont\def\encodingdefault{LGR}}
\DeclareRobustCommand{\textgreek}[1]{\leavevmode{\greektext #1}}
\ProvideTextCommand{\~}{LGR}[1]{\char126#1}

\makeatother

\usepackage{babel}
\begin{document}
\title{Simple formulas of $\pi$ in terms of $\Phi$}
\author{Angelo Pignatelli}
\date{April 9th, 2024}
\maketitle
\begin{abstract}
A re-calculation of a known family of formulas of $\pi$ is carried
out, revisiting the old Archimedes' algorithm. This allows to identify
a general family equation and three new simple formulas of $\pi$
in terms of the golden ratio $\Phi$ in the form of infinite nested
square roots, with some geometrical properties that enhance the link
between the circle and the golden ratio. Applying the same criteria,
a fourth formula is given, that brings to the known Dixon's squaring
the circle approximation, thus an easier approach to this problem
is suggested, by a rectangle with both sides proportional to the golden
ratio $\Phi$.
\end{abstract}

\section*{Introduction}

In last century it was challenging and interesting to find formulas
of $\pi$ in terms of the golden ratio $\Phi$ (so involving together
two of most famous irrational constants), without transcendental functions,
as the well known $\pi=\frac{10}{3}\arcsin\left(\frac{\phi}{2}\right)$.
Among these works, few examples are here reported:

\begin{equation}
\pi=\frac{5\sqrt{2+\phi}}{2\phi}\sum_{n=0}^{\infty}\left(\frac{1}{2\phi}\right)^{5n}\left(\frac{1}{5n+1}+\frac{1}{2\phi^{2}\left(5n+2\right)}-\frac{1}{2^{2}\phi^{3}(5n+3)}-\frac{1}{2^{3}\phi^{3}(5n+4)}\right)
\end{equation}
Equation (1) has been presented by Chan in {[}1{]}, inspired by the
work of Bailey, Borwein and Plouffe (so called BBP-formulas) in {[}2{]}.
\begin{equation}
\frac{\pi^{2}}{50}=\sum_{k=0}^{\infty}\left(\frac{\phi^{2}}{\left(5k+1\right)^{2}}-\frac{\phi}{\left(5k+2\right)^{2}}-\frac{\phi^{2}}{\left(5k+3\right)^{2}}+\frac{\phi^{5}}{\left(5k+4\right)^{2}}+\frac{2\phi^{2}}{\left(5k+5\right)^{2}}\right)\phi^{-5k}
\end{equation}
Equation (2) has been discovered by B. Cloitre and reported by Chan
in {[}3{]}, also inspired by BBP-formulas. 

The aim of this work, presented in next pages, is focused on identifying
other simpler formulas of $\pi$ in terms of $\Phi$.

\begin{figure}
\begin{centering}
\includegraphics[scale=0.35]{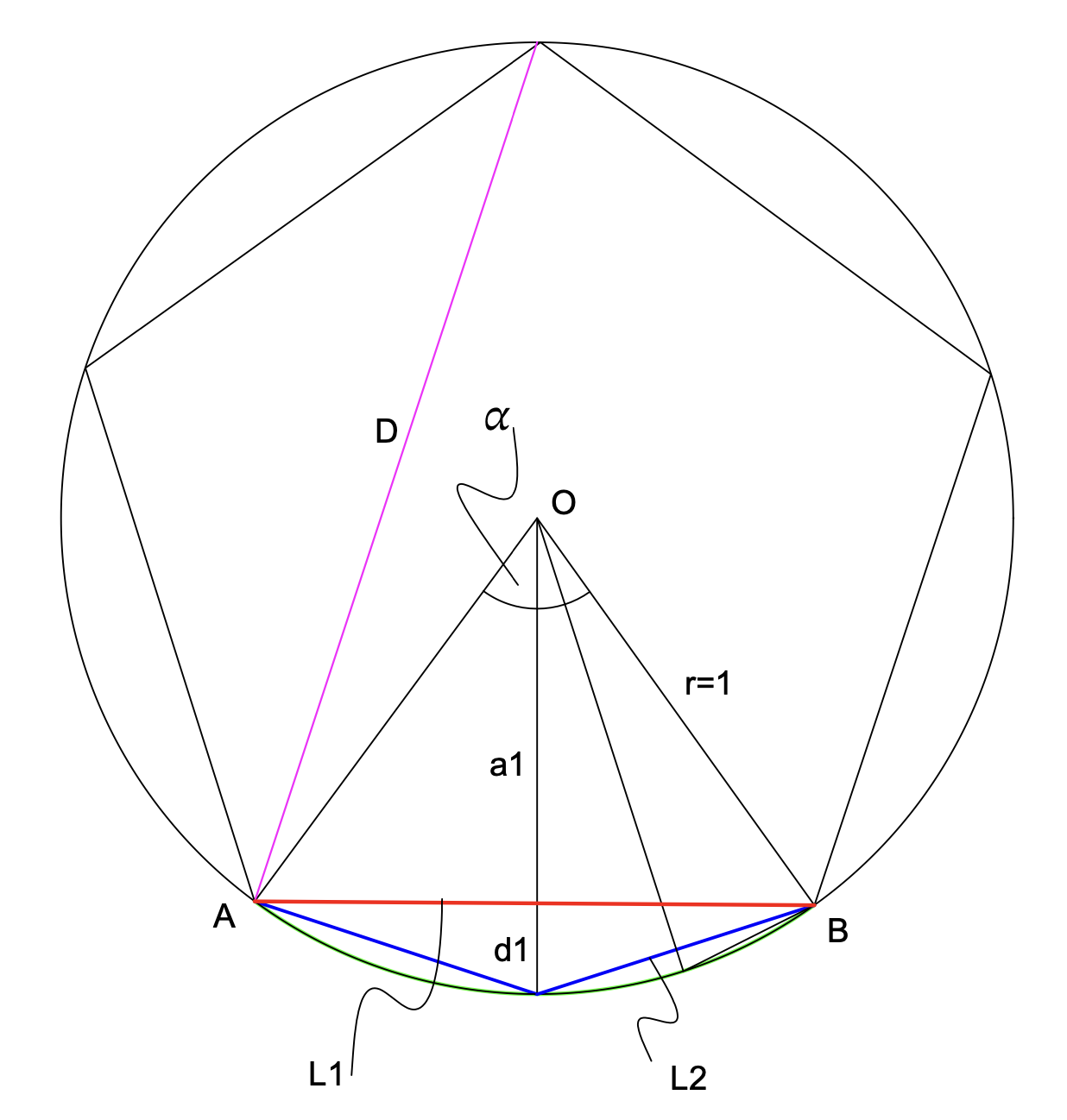}
\par\end{centering}
\caption{Iterative halfway bisection of chord and arc - example on a pentagon}

\end{figure}

\section*{Nested square roots formulas of $\pi$}

In order to show other simple formulas of $\pi$ in terms of $\Phi$,
first it needs to easily share the calculation behind the family of
the known formulas of $\pi$ in the form of nested square roots. The
approach starts from the idea of Archimedes, resumed several centuries
later by F. Viete (as reported by Beckmann in {[}4{]}), then more
recently completely restructured by Servi in {[}5{]}.

Let us start considering a regular polygon inscribed in a circle with
unitarian radius, with N sides ($N\mathbb{\epsilon N},N\geq3$) of
length L$_{1}$$(L_{1}=2\sin\frac{\alpha}{2}$, $\alpha=\frac{2\pi}{N}$
), with perimeter $P_{1}=NL_{1}$, as in Figure 1 (where, as example,
is represented a pentagon).

With the polygon obtained doubling the sides, considering $a_{1}=\sqrt{1-\left(\frac{L_{1}}{2}\right)^{2}}$,
$d_{1}=1-a_{1}$, the perimeter become $P_{2}=N2L_{2}$ with
\begin{equation}
L_{2}=\sqrt{\frac{L_{1}^{2}}{4}+\left(1-\sqrt{1-\frac{L_{1}^{2}}{4}}\right)^{2}}=\sqrt{2-2\sqrt{1-\frac{L_{1}^{2}}{4}}}=\sqrt{2-\sqrt{4-L_{1}^{2}}}
\end{equation}

Re-iterating the process, doubling the polygon at each step, from
Equation (3) the following succession is obtained
\begin{equation}
\begin{cases}
L_{1}=2\sin\frac{\pi}{N} & N\epsilon\mathbb{N},N\geq3\\
L_{n}=\sqrt{2-\sqrt{4-L_{n-1}^{2}}} & n\epsilon\mathbb{N},n\geq2\\
\pi=N\lim_{n\rightarrow\infty}2^{n-2}L_{n}
\end{cases}
\end{equation}

Expanding Equations (3) and (4):

\[
L_{3}=\sqrt{2-\sqrt{4-L_{2}^{2}}}=\sqrt{2-\sqrt{4-\left(2-\sqrt{4-L_{1}^{2}}\right)}}=\sqrt{2-\sqrt{2+\sqrt{4-L_{1}^{2}}}}
\]

\[
L_{4}=\sqrt{2-\sqrt{4-L_{3}^{2}}}=\sqrt{2-\sqrt{4-\left(2-\sqrt{2+\sqrt{4-L_{1}^{2}}}\right)}}=\sqrt{2-\sqrt{2+\sqrt{2+\sqrt{4-L_{1}^{2}}}}}
\]

\begin{equation}
L_{n}=\sqrt{2-\sqrt{2+\sqrt{2+\sqrt{2+\ldots\sqrt{2+\sqrt{4-L_{1}^{2}}}}}}}\frac{}{}\frac{}{}\frac{}{}n\frac{}{}square\frac{}{}roots
\end{equation}

Using the results of Equation (5), the known family of equations of
\ensuremath{\pi} in the form of continued square roots follows, valid
for any regular polygons with N sides and side length $L_{1}$.

\begin{equation}
\begin{cases}
L_{1}=2\sin\frac{\pi}{N} & N\epsilon\mathbb{N},N\geq3\\
\pi=N\lim_{n\rightarrow\infty}2^{n-1}\sqrt{2-\sqrt{2+\sqrt{2+\sqrt{2+\ldots\sqrt{2+\sqrt{4-L_{1}^{2}}}}}}} & n+1\frac{}{}square\frac{}{}roots
\end{cases}
\end{equation}

Now we would like to extend this approach, freeing from the regular
polygons.

Let us focus on the arc $\overset{}{\overset{\smallfrown}{AB}=\alpha}<\pi\frac{}{}rad$
, implied with the chord $\overline{AB}=2\sin\frac{\alpha}{2}$ (Figure
1); with this in mind, we can call $\mu\epsilon\mathbb{R}$ the ratio
between the circumference length and the arc $\overset{\smallfrown}{AB}$:
\begin{equation}
\mu=\frac{2\pi}{\alpha}\frac{}{}\frac{}{}\rightarrowtail\frac{}{}\frac{}{}\pi=\frac{1}{2}\mu\alpha\frac{}{}\frac{}{}\frac{}{}\frac{}{}\frac{}{}\frac{}{}with\frac{}{}0<\alpha<\pi\frac{}{}\frac{}{}\frac{}{}\rightarrowtail\frac{}{}\frac{}{}2<\mu<+\infty
\end{equation}

Applying a similar strategy, it is proven that, dividing iteratively
halfway the arc $\overset{\smallfrown}{AB}$ (as in Figure 1), the
sum of the chords implied converge to the length of the arc $\overset{\smallfrown}{AB}\overset{}{=\alpha}$
, thus
\begin{equation}
\begin{cases}
\mu=\frac{2\pi}{\alpha} & 0<\alpha<\pi\\
\pi=\mu\lim_{n\rightarrow\infty}2^{n}\sin\left(\frac{\alpha}{2^{n}}\right)
\end{cases}
\end{equation}

Applying several times the goniometric bisection formulas $\sin\frac{\gamma}{2}=\sqrt{\frac{1-\cos\gamma}{2}}$
, $\cos\frac{\gamma}{2}=\sqrt{\frac{1+\cos\gamma}{2}}$, and ones
the formula $\sin^{2}\gamma+\cos^{2}\gamma=1$ to the expression $\sin\frac{\alpha}{2^{n}}$
we obtain:
\begin{flushleft}
\[
\sin\left(\frac{\alpha}{2^{n}}\right)=\sin\left(\frac{\frac{\alpha}{2^{n-1}}}{2}\right)=\sqrt{\frac{1-\cos\left(\frac{\alpha}{2^{n-1}}\right)}{2}}=\frac{1}{2}\sqrt{2-2\cos\left(\frac{\alpha}{2^{n-1}}\right)}=\frac{1}{2}\sqrt{2-2\cos\left(\frac{\frac{\alpha}{2^{n-2}}}{2}\right)}=
\]
\[
\frac{}{}\frac{}{}\frac{}{}\frac{}{}\frac{}{}=\frac{1}{2}\sqrt{2-2\sqrt{\frac{1+\cos\left(\frac{\alpha}{2^{n-2}}\right)}{2}}}=\frac{1}{2}\sqrt{2-\sqrt{2+2\cos\left(\frac{\alpha}{2^{n-2}}\right)}}=
\]
\[
=\frac{1}{2}\sqrt{2-\sqrt{2+\sqrt{2+\ldots\sqrt{2+2cos\left(\frac{\alpha}{2}\right)}}}}=\frac{1}{2}\sqrt{2-\sqrt{2+\sqrt{2+\ldots\sqrt{2+2\sqrt{1-\sin^{2}\left(\frac{\alpha}{2}\right)}}}}}=
\]
\[
=\frac{1}{2}\sqrt{2-\sqrt{2+\sqrt{2+\ldots\sqrt{2+\sqrt{4-\left[2\sin\left(\frac{\alpha}{2}\right)\right]^{2}}}}}}
\]
\par\end{flushleft}

Using this result in Equation (8), noting that $L_{1}=\overline{AB}=2\sin\left(\frac{\alpha}{2}\right)$,
finally the following general formula is obtained:

\begin{equation}
\begin{cases}
\mu=\frac{2\pi}{\alpha} & 0<\alpha<\pi\frac{}{}\frac{}{},\frac{}{}\frac{}{}\mu\epsilon\mathbb{R},\mu>2\\
L_{1}=2\sin\left(\frac{\alpha}{2}\right)\\
\pi=\mu\lim_{n\rightarrow\infty}2^{n-1}\sqrt{2-\sqrt{2+\sqrt{2+\ldots\sqrt{2+\sqrt{4-L_{1}^{2}}}}}} & n+1\frac{}{}square\frac{}{}roots
\end{cases}
\end{equation}

Equation (9) generalizes and confirms Equation (6). Equation (9) coincides
with Equation (6) when $\mu=N\epsilon\mathbb{N}$, or in other words
when \ensuremath{\alpha} is an integer divisor of the circle, and
in this case $L_{1}=2\sin\left(\frac{\alpha}{2}\right)=2\sin\left(\frac{\pi}{N}\right)$.
Also the succession in Equation (4) can be extended substituting N
with $\mu=\frac{2\pi}{\alpha}$ and $L_{1}$ with $L_{1}=2\sin\left(\frac{\alpha}{2}\right)$
for values not integer of $\frac{2\pi}{\alpha},$ with $\alpha\epsilon\left(0,\pi\right)$.

Equation (9), as Equation (6), arises some interest when applied for
values not transcendental of $\sin\left(\frac{\alpha}{2}\right)$.
Some of these instances follow.

Applying Equation (6) to an equilateral triangle ($N=3,L_{1}=2\sin\frac{\pi}{3}=\sqrt{3}$)
or an hexagon ($N=6,L_{1}=2\sin\frac{\pi}{6}=1$), or a dodecagon
($N=12,L_{1}=2\sin\frac{\pi}{12}=\frac{\sqrt{6}-\sqrt{2}}{2}$ ) the
following Equation (10), coinciding with the formula (3) presented
by Servi in {[}5{]}, is obtained:

\begin{equation}
\pi=3\lim_{n\rightarrow\infty}2^{n-1}\sqrt{2-\sqrt{2+\sqrt{2+\ldots\sqrt{2+\sqrt{3}}}}}\frac{}{}\frac{}{},n\frac{}{}square\frac{}{}roots
\end{equation}

Applying Equation (6) to a square (N=4, $L_{1}=2\sin\frac{\pi}{4}=\sqrt{2}$),
it is possible to obtain the following Equation (11), coinciding with
the formula (1) presented by Servi in {[}5{]}:

\begin{equation}
\pi=\lim_{n\rightarrow\infty}2^{n+1}\sqrt{2-\sqrt{2+\sqrt{2+\ldots\sqrt{2+\sqrt{2}}}}}\frac{}{}\frac{}{},n+1\frac{}{}square\frac{}{}roots
\end{equation}

Applying Equation (9) to the arc $\overset{\smallfrown}{AB}=\alpha=\frac{3}{4}\pi$
(135°), $\mu=\frac{2\pi}{\alpha}=\frac{8}{3}$, calculating $L_{1}=2\sin\left(\frac{\alpha}{2}\right)=\sqrt{2+\sqrt{2}}$
and the last square root in Equation (9) $\sqrt{4-L_{1}^{2}}=\sqrt{4-(2+\sqrt{2)}}=\sqrt{2-\sqrt{2}},$
the following Equation (12), coinciding with the formula (2) presented
by Servi in {[}5{]}, is obtained:
\begin{equation}
\pi=\frac{1}{3}\lim_{n\rightarrow\infty}2^{n+2}\sqrt{2-\sqrt{2+\sqrt{2+\ldots\sqrt{2+\sqrt{2-\sqrt{2}}}}}}\frac{}{}\frac{}{}\frac{}{}n+2\frac{}{}square\frac{}{}roots
\end{equation}

Applying Equation (9) to the arc $\overset{\smallfrown}{AB}=\alpha=\frac{5}{6}\pi$
(150°), $\mu=\frac{2\pi}{\alpha}=\frac{12}{5}$, calculating $L_{1}=2\sin\left(\frac{\alpha}{2}\right)=2\frac{\sqrt{6}+\sqrt{2}}{4}=\frac{\sqrt{6}+\sqrt{2}}{2}$
and the last square root in Equation (9) $\sqrt{4-L_{1}^{2}}=\sqrt{4-\frac{\left(\sqrt{6}+\sqrt{2}\right)^{2}}{4}}=\sqrt{4-\frac{6+2+2\sqrt{12}}{4}}=\frac{1}{2}\sqrt{8-4\sqrt{3}}=\sqrt{2-\sqrt{3}},$
the following Equation (12), coinciding with the formula (4) presented
by Servi in {[}5{]}, is obtained:
\begin{equation}
\pi=\frac{3}{5}\lim_{n\rightarrow\infty}2^{n+1}\sqrt{2-\sqrt{2+\sqrt{2+\ldots\sqrt{2+\sqrt{2-\sqrt{3}}}}}}\frac{}{}\frac{}{}\frac{}{}n+2\frac{}{}square\frac{}{}roots
\end{equation}

\section*{New simple formulas of $\pi$ in terms of $\phi$}

Now applying specifically Equations (6) and (9), the following three
simple formulas of $\pi$ in terms of $\Phi$ are identified:

\begin{equation}
\pi=\frac{5}{2}\lim_{n\rightarrow\infty}2^{n}\sqrt{2-\sqrt{2+\sqrt{2+\ldots\sqrt{2+\sqrt{2+\phi}}}}}\frac{}{}\frac{}{}\frac{}{}\frac{}{}\frac{}{}\frac{}{}\frac{}{}\frac{}{}\frac{}{}n\frac{}{}square\frac{}{}roots
\end{equation}

\begin{equation}
\pi=\frac{5}{3}\lim_{n\rightarrow\infty}2^{n}\sqrt{2-\sqrt{2+\sqrt{2+\ldots\sqrt{2+\sqrt{3-\phi}}}}}\frac{}{}\frac{}{}\frac{}{}n+1\frac{}{}square\frac{}{}roots
\end{equation}

\begin{equation}
\pi=\frac{5}{4}\lim_{n\rightarrow\infty}2^{n}\sqrt{2-\sqrt{2+\sqrt{2+\ldots\sqrt{2+\sqrt{2-\phi}}}}}\frac{}{}\frac{}{}\frac{}{}n+1\frac{}{}square\frac{}{}roots
\end{equation}

Equation (14) is obtained applying the Equation (6) to a pentagon
inscribed in the circumference (as in Figure 1) with unitarian radius
($N=5$, $\alpha=\frac{2}{5}\pi$, with simple passages $L_{1}=2\sin\left(\frac{\pi}{5}\right)=2\frac{\sqrt{10-2\sqrt{5}}}{4}=\sqrt{\frac{10-2\sqrt{5}}{4}}=\sqrt{\frac{5-\sqrt{5}}{2}}=\sqrt{3-\frac{1+\sqrt{5}}{2}}=\sqrt{3-\phi}$).
From Equation (6), and remembering that $\phi^{n}=\phi^{n-1}+\phi^{n-2}$
, the last two square roots become $\sqrt{2+\sqrt{4-L_{1}^{2}}}=\sqrt{2+\sqrt{4-\left(3-\phi\right)}}=\sqrt{2+\sqrt{1+\phi}}=\sqrt{2+\phi}$.
Since between the diagonal D and the side L of a pentagon results
$D=\phi L$ (Ghyka in {[}6{]}), it is noted the fact quite singular
that in this case the last square root in Equation (14) represents
exactly the length of the diagonal $D_{1}=\phi L_{1}=\phi\sqrt{3-\phi}=\sqrt{3\phi^{2}-\phi^{3}=}\sqrt{3\phi^{2}-\left(\phi^{2}+\phi\right)}=\sqrt{2\phi^{2}-\phi}=\sqrt{2\left(\phi+1\right)-\phi}=\sqrt{2+\phi}$.

Equation (15) is obtained applying the Equation (9) to the arc $\overset{\smallfrown}{AB}=\alpha=\frac{3}{5}\pi$
(108°), $\mu=\frac{2\pi}{\alpha}=\frac{5}{2}$. Calculating $L_{1}=2\sin\left(\frac{\alpha}{2}\right)=2\frac{1+\sqrt{5}}{4}=\frac{1+\sqrt{5}}{2}=\phi$
(in fact in this case the chord $\overline{AB}$ with two radius forms
the triangle gnomon of the golden triangle), thus the last square
root in Equation (9) become $\sqrt{4-L_{1}^{2}}=\sqrt{4-\phi^{2}}=\sqrt{4-\left(\phi+1\right)}=\sqrt{3-\phi}$.
It is noted the fact quite singular that in this case the last square
root in the formula Equation (16) represents exactly the length of
the side of a pentagon inscribed in the circle.

Equation (16) is obtained applying the Equation (9) to the arc $\overset{\smallfrown}{AB}=\alpha=\frac{4}{5}\pi$
(144°), $\mu=\frac{2\pi}{\alpha}=\frac{10}{3}$. Calculating $L_{1}=2\sin\left(\frac{\alpha}{2}\right)=\frac{1}{2}\sqrt{10+2\sqrt{5}}=\sqrt{\frac{4+1+\sqrt{5}}{2}}=\sqrt{2+\phi}$
(coinciding with the diagonal of a pentagon inscribed in the circle),
thus the last square root in Equation (9) become $\sqrt{4-L_{1}^{2}}=\sqrt{2-\phi}$.
It is noted the fact quite singular that in this case the last square
root in Equation (16) represents exactly the length of the side of
a decagon inscribed in the circle (as $2\sin\frac{36{^\circ}}{2}=\frac{1}{2}\left(\sqrt{5}-1\right)=\sqrt{\left(\frac{\sqrt{5}-1}{2}\right)^{2}}=\sqrt{\frac{6-2\sqrt{5}}{4}}=\sqrt{\frac{4-1-\sqrt{5}}{2}}=\sqrt{2-\phi}$).

\section*{Geometrical properties}

\begin{figure}
\begin{centering}
\includegraphics[scale=0.35]{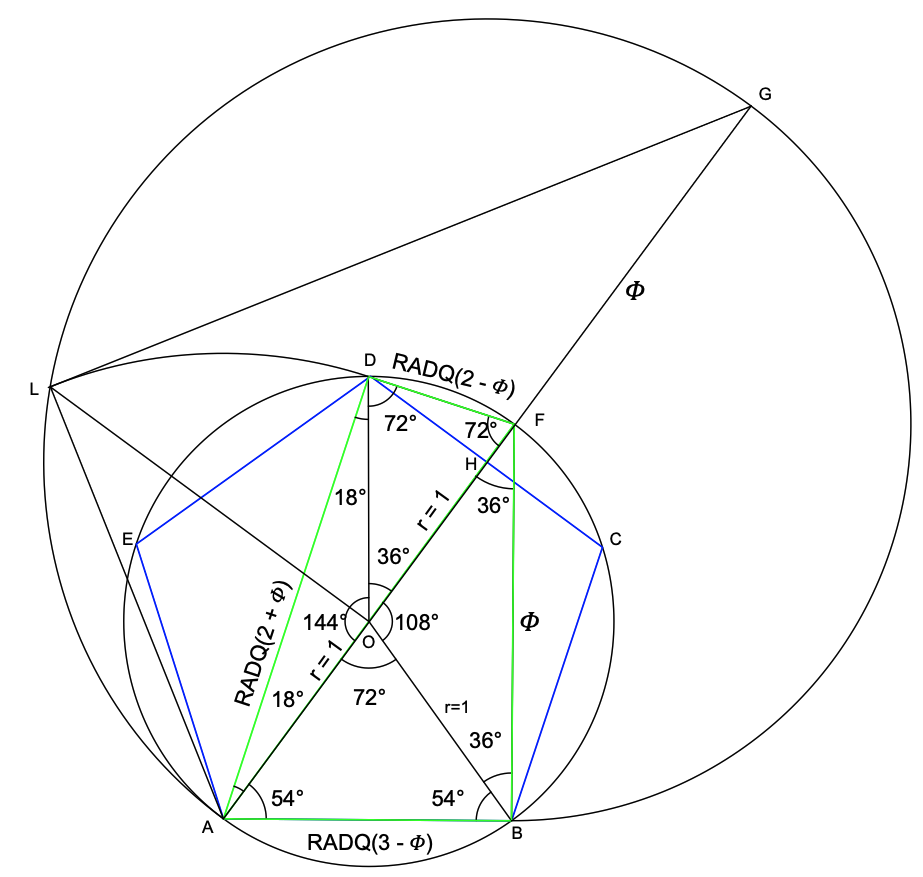}
\par\end{centering}
\caption{Geometrical properties between the circle, \ensuremath{\pi} and \ensuremath{\phi}}

\end{figure}

Figure 2 helps to understands the singularities pointed out in last
square root of Equations (14), (15) and (16).

For Equation (14), as $L_{1}=\overline{AB}=\sqrt{3-\phi}$ (side of
the pentagon inscribed in a circle with unitarian radius), applying
the last square root $\sqrt{4-L_{1}^{2}}$ means finding the other
cathetus $\overline{BF}$ of the right triangle ABF, being the hypotenuse
$\overline{AF}=2$, then $\overline{BF}=\sqrt{4-(3-\Phi)}=\sqrt{1+\varPhi}=\Phi$
(also as side of a gnomon triangle OBF of the golden triangle). Then,
for the following square root $\sqrt{2+\phi}$, with the graphical
approach it is possible to construct the segment AG with length $2+\phi$,
then drawing half a circle with center in H e diameter AG, $\overline{AL}$
results the square root of $\overline{AG}$ (from the equivalence
of the right triangles AGL and AOL, with $\overline{AO}=1$); finally
we get the diagonal $\overline{AD}=\overline{AL}=\sqrt{2+\varPhi}$.

For Equation (15), as $L_{1}=\overline{BF}=\phi$ (side of the gnomon
triangle), applying the last square root $\sqrt{4-L_{1}^{2}}$ means
finding the other cathetus $\overline{AB}$ of the right triangle
ABF, of length $\sqrt{3-\phi}$ (side of the penthagon).

For Equation (16), as $L_{1}=\overline{AD}=\sqrt{2+\phi}$ (diagonal
of the pentagon), applying the last square root $\sqrt{4-L_{1}^{2}}$
means finding the other cathetus $\overline{DF}$ of the right triangle
AFD, then $\overline{DF}=\sqrt{4-(2+\Phi)}=\sqrt{2-\Phi}$ (side of
the decagon).

From this geometrical approach it is evident that, if we apply Equation
(9) with $L_{1}=\overline{DF}=\sqrt{2-\varPhi}$ (side of a decagon),
applying the last square root $\sqrt{4-L_{1}^{2}}$ means finding
the other cathetus $\overline{AD}$ of the right triangle AFD, of
length $\sqrt{2+\phi}$ (diagonal of the pentagon), thus with the
same result obtained in Equation (14) starting with $L_{1}=\overline{AB}=\sqrt{3-\varPhi}$
(side of a pentagon), as we could expect by doubling the sides of
the pentagon on the first iteration of Equation (4).

Referring to Figure 2, it should be noted that $\overline{AB},$ $\overline{BF},$
$\overline{AD}$ are diagonals of the decagon inscribed in the circle,
with side $\overline{DF}.$

Equations (14), (15), (16) and geometrical properties in Figure 2,
identified in a circle with unitarian radius and its inscribed pentagon
(and decagon), arise some relationships between the circle, $\pi$
and $\Phi$, extending the relations inside the pentagon and its pentagram
constructed with its diagonals (in Figure 3), well known since Pythagoras
ancient times (Ghyka in {[}6{]}).

\begin{figure}
\begin{centering}
\includegraphics[scale=0.36]{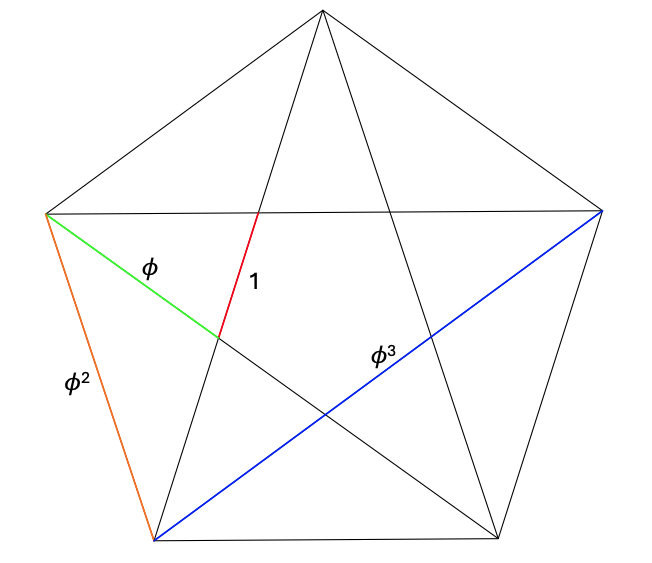}
\par\end{centering}
\caption{Relations in Pythagorean pentagram}

\end{figure}

\section*{New approximate ``squaring'' the circle proposal}

Let us consider a particular angle, the arc $\overset{\smallfrown}{AB}=\alpha=\frac{2\pi}{\phi^{2}}$
(\ensuremath{\eqsim} 137.5° called also ``golden angle'', that is
found many times in nature, for instance in phyllotaxis), with $\mu=\frac{2\pi}{\alpha}=\phi^{2}$.

Calculating $L_{1}=2\sin\left(\frac{\alpha}{2}\right)=2\sin\left(\frac{\pi}{\phi^{2}}\right)\simeq1.86406485$,
Equation (9) becomes:
\begin{equation}
\pi=\phi^{2}\lim_{n\rightarrow\infty}2^{n-1}\sqrt{2-\sqrt{2+\sqrt{2+\ldots\sqrt{2+2\sqrt{1-L_{1}^{2}}}}}}\frac{}{}\frac{}{}\frac{}{}n+1\frac{}{}square\frac{}{}roots
\end{equation}

The interest on this formula arises noticing that the limit converges
to a number 1.199981546..., thus can be approximate to 1.2 with an
error lower than 0.00005.

The approximation of Equation (17) provides a mathematical source
to the well known approximate formula in following Equation (18) between
\ensuremath{\pi} and \ensuremath{\phi}:
\begin{equation}
\pi\simeq\frac{6}{5}\phi^{2}=\frac{6}{5}\left(1+\Phi\right)=3,141640...\frac{}{}\frac{}{}\frac{}{}\frac{}{}\frac{}{}\frac{}{}err<0.00005
\end{equation}

This approximation has been pointed out by Dixon in {[}7{]}, also
showing an interesting procedure to draw a square with an area of
$\frac{6}{5}\left(1+\phi\right)\simeq\pi$ with an error lower than
0.00005 (in Figure 4). It could be interesting to mention that the
relation in Equation (18) is known at least from the 12th century
by the French master masons that built the gothic cathedrals, as proved
by Frederic in {[}8{]}.
\begin{figure}
\begin{centering}
\includegraphics[scale=0.35]{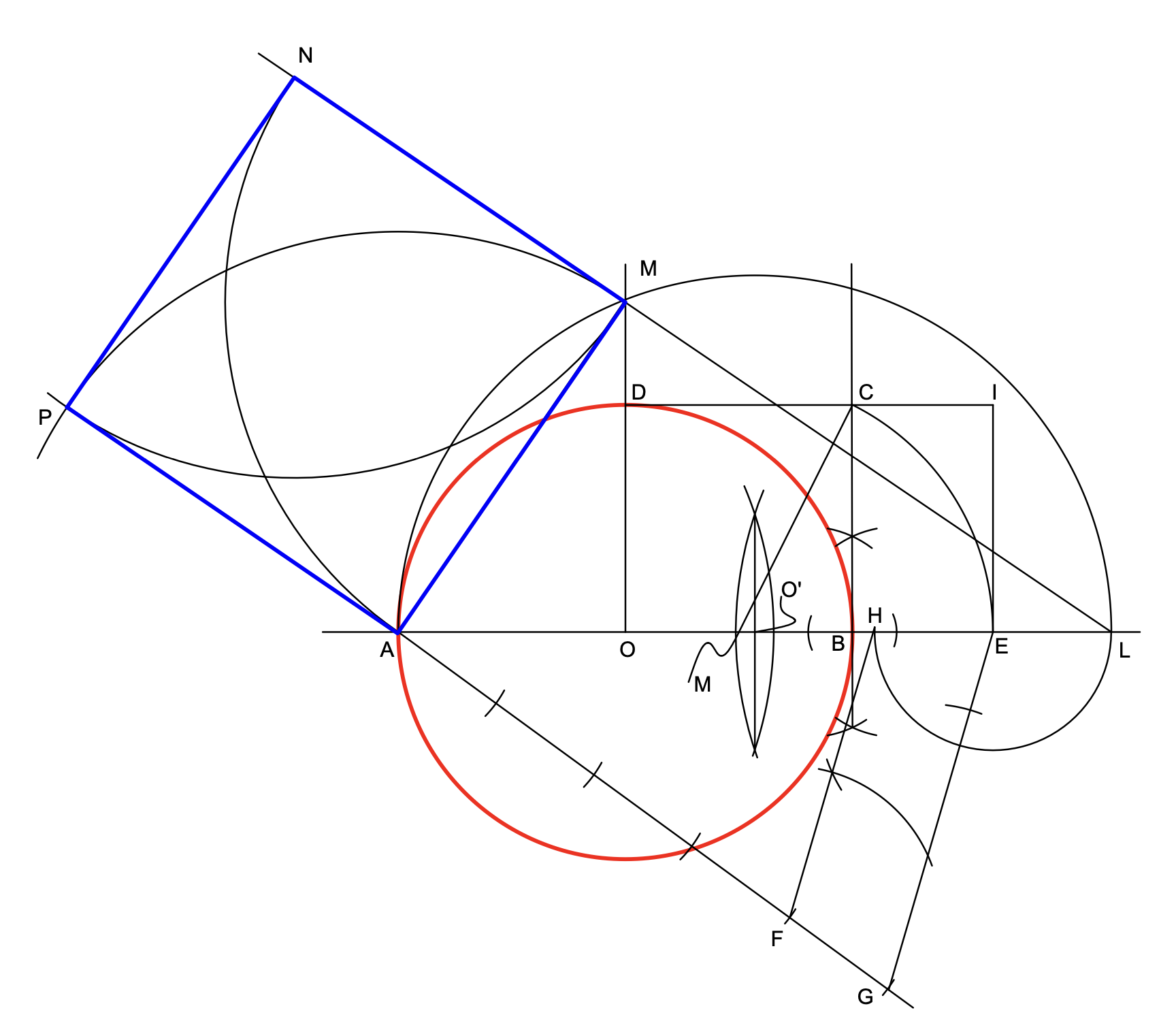}
\par\end{centering}
\caption{Approximate squaring the circle by R. A. Dixon}

\end{figure}

\begin{quote}
Procedure for the constructions just by rule and compass (Figure 4):
1. draw a circle with radius = 1; 2. trace the golden rectangle OEID;
3. apply the rule to divide a segment in five equal parts with segment
AG, identifying the fifth part HE of $\overline{AE}=\overline{AO}+\overline{OE}=1+\phi$,
then add this 1/5 to the right, in order to identify the segment AL
with length $\frac{6}{5}(1+\phi)$; 4. trace half a circle on the
diameter AL finding point M as the intersection with the vertical
line from the centre O; 5. construct the square AMNP on the segment
AM. As the triangles ALM and AOM are similar, $\overline{AL}:\overline{AM}=\overline{AM}:\overline{AO}\rightarrow\overline{AM}^{2}=\overline{AL}\rightarrow\overline{AM}=\sqrt{\overline{AL}}=\sqrt{\frac{6}{5}\left(1+\phi\right)}\simeq\sqrt{\pi}$.
\end{quote}
We propose here another easier way to approximate the ``squaring''
the circle based on Equation (18) with not a square but a rectangle,
with sides length $\phi$ and $\frac{6}{5}\phi$, whose area $\frac{6}{5}\phi^{2}$
is quite close (with error lower than 0.00005) to the area $\pi$
of the circle, as in Figure 5.
\begin{figure}
\begin{centering}
\includegraphics[scale=0.35]{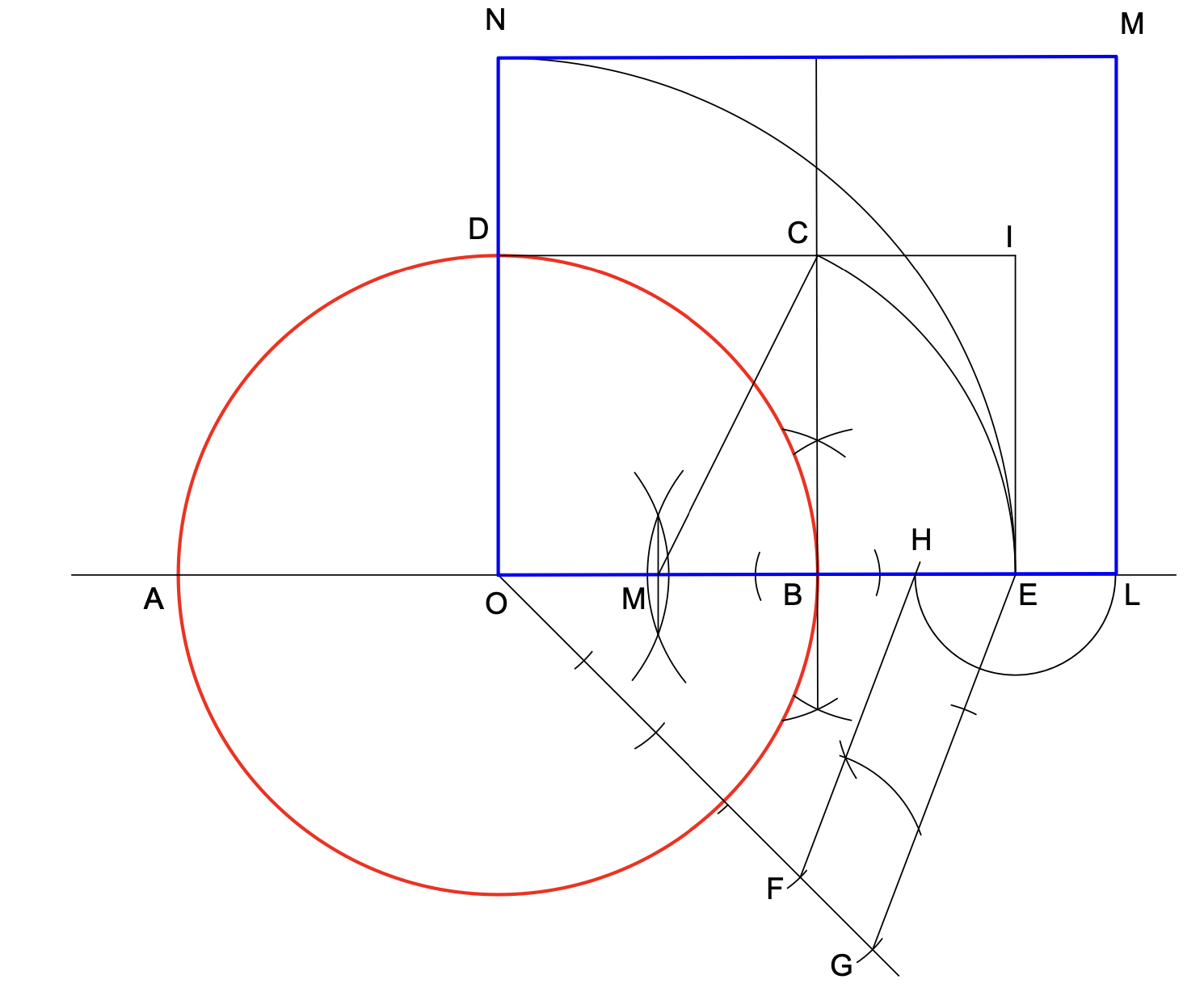}
\par\end{centering}
\caption{Approximate ``squaring'' the circle by a ``$\pi\sim$ rectangle''}
\end{figure}

\begin{quote}
Procedure for the constructions just by rule and compass (Figure 5):
1-2. apply the same previous steps; 3. apply the rule to divide a
segment in five equal parts with segment OG, identifying the fifth
part HE of $\overline{OE}=\phi$, then add this 1/5 to the right,
in order to identify the segment OL with length $\frac{6}{5}\phi$;
4. trace the arc with center on O from E in order to identify the
point N as intersection with the vertical line from O, so $\overline{ON}=\phi$;
5. Trace the rectangle OLMN, (that we can call the ``$\pi\backsim$rectangle''),
that has sides length $\phi$ x $\frac{6}{5}\phi$ , thus with area
$\frac{6}{5}\phi^{2}\simeq\pi$.
\end{quote}

\section*{Conclusions}

After sharing the calculation behind the family of the known formulas
of $\pi$ in the form of nested square roots, with the presented general
formula in Equation (9) three new simple formulas of $\pi$ in terms
of $\Phi$ are given in Equations (14), (15) and (16), arising some
interesting geometrical properties (in Figure 2), that enhance the
link between the circle and the golden ratio; these relationships
could be deeper investigated in the future. A mathematical basis,
in Equation (17), is provided for the well known approximation of
$\pi$ in terms of $\Phi$ in Equation (18); a so called ``$\pi\backsim$
rectangle'', that has sides length $\phi$ x $\frac{6}{5}\phi$ ,
is suggested as an approximated ``squaring'' the circle problem
(Figure 5).

\section*{References}
\begin{description}
\item [{{[}1{]}}] Chan, C. (2004) \textgreek{p} in terms of \ensuremath{\phi},
Mathematical Sciences Program, University of Illinois at Springfield,
Springfield, IL 62703-5407.
\item [{{[}2{]}}] Bailey, D. H., Borwein, P. B. and Plouffe, S. (1997)
On the Rapid Computation of Various Polylogarithmic Constants'',
Math. Comp., 66, 903--913.
\item [{{[}3{]}}] Chan, C. (2009) Machin-type formulas expressiong \textgreek{p}
in terms of phi, The Fibonacci Quarterly 46/47, 32-37.
\item [{{[}4{]}}] Beckmann, (1989) A History of Pi, Dorset Press, 3rd ed.,
New York.
\item [{{[}5{]}}] Servi, D. (2003) Nested square roots of 2, Amer. Math.
Monthly 110 (4) 326-330.
\item [{{[}6{]}}] Ghyka, C. (2016) The Golden Number, Inner Tradition,
Rochester.
\item [{{[}7{]}}] Dixon, (1991) Squaring the circle, Mathographics Blackwell
1987, pp. 44--47, reprinted by Dover Publications.
\item [{{[}8{]}}] Frédéric, (2022) Squaring the circle like a medieval
master mason, UNSW School of Mathematics and Statistics, ParabolaVolume
58 Issue 2.
\end{description}
\begin{doublespace}
Author:
\end{doublespace}

Ph.D. Angelo Pignatelli

email: ing@angelopignatelli.it
\end{document}